\documentclass[a4paper]{article}

\usepackage{graphicx,epstopdf}
\usepackage{amsfonts,amsmath,amssymb,amsthm}
\usepackage{graphicx}
\usepackage{mathrsfs}
\usepackage[title, titletoc, header]{appendix}

\newtheorem{theorem}{Theorem}[section]

\usepackage{amsmath}

\usepackage{graphicx,epstopdf}

\newcommand{\NN}{\mathbb{N}}
\newcommand{\RR}{\mathbb{R}}

\begin{document}

\begin{center}
\large{ \textbf{ {\Large Reliable approximation of separatrix manifolds 
\vskip 0.1cm
in  competition models with safety niches}}}
\end{center}

\begin{center}
Roberto Cavoretto, Alessandra De Rossi, Emma Perracchione and \\ Ezio Venturino
\end{center}

\begin{center}
Department of Mathematics "G. Peano", University of Turin - Italy
\end{center}
\vskip 0.5cm

\textbf{Abstract}
In dynamical systems saddle points partition the domain into basins of attractions of the remaining locally stable equilibria. This situation
is rather common especially in population dynamics models, like prey-predator or competition systems. Focusing on squirrels
population models with niche, in this paper we design algorithms for the detection and the refinement of points lying on the separatrix manifold 
partitioning the phase space. We consider both the two populations and the three populations cases.
To reconstruct the separatrix curve and surface,
we apply the Partition of Unity method, which makes use of Wendland's functions as local approximants.


\section{Introduction}

The competition of two or more species that live in the same environment can be modelled mathematically
by a differential system, whose unknowns represent the populations sizes as functions of time. Their interactions are then described by a number of parameters (see \cite{Arrowsmith90,Brauer01,Murray02,Wiggins03}).
To obtain a particular solution of the system, we need to know the initial state of the system.
The system then in general evolves toward stable equilibria, for suitable parameter choices. In classical two population competition models
indeed limit cycles are excluded and further the principle of competitive exclusion applies, for which the ``best fit'' population
outcompetes the other one. The knowledge of the winner depends on the system's initial conditions. If the latter lie in what is called the
basin of attraction of a certain equilibrium point, the final population configuration will be the one at this specific equilibrium.
Therefore, it is important to assess the basins of attraction of each possible equilibrium.
We can thus imagine to partition the phase space in different regions depending on where the trajectory originating in
them will ultimately stabilize.
The aim of this work is to construct an approximation curve or an approximation surface,
which partitions the considered domain in two or more regions, called the basins of attraction of each equilibrium.

In particular, in this article we discuss two specific population models with niche, which investigate squirrels competition of ecosystems
composed of two and three different populations (see also \cite{derossi13,grosso12,gurnell04}).
The former considers a two population model with competition between red native and grey exotic squirrels,
while the latter involves a three population model with competition among red native, red indigenous and grey exotic squirrels.
At first, we carry out an analytical study of the two models, aimed at finding the location of equilibrium points.
We establish conditions to be imposed on the parameters so that the behavior described above in fact occurs and the separatrix manifolds exist.
Then, after choosing parameters which satisfy these assumptions for feasibility and stability of the equilibria,
we proceed to approximate the separatrix curve and surface. For this purpose we have implemented several \textsc{Matlab}
functions for the approximation of the points, obtained by a bisection algorithm, and the graphical representation of the separating curve and surface.
After detecting the points lying on the separatrix manifold, we first apply a refinement algorithm in order to reduce the number of points to interpolate.
We then approximate the curve and surface using the Partition of Unity method using as local approximants the compactly supported Wendland's functions
(see, e.g., \cite{Fasshauer07,Wendland05}). The latter is an effective and efficient tool in approximation theory,
since it allows to interpolate a large number of scattered data in an accurate and stable way (see \cite{Allasia11,Cavoretto12a,cavoretto12c,cavoretto12d}).
For a two populations model, \cite{cavoretto11}, a different refinement algorithm has already been used.

The paper is organized as follows. In Sections 2 and 3 we consider the two populations and three populations models respectively,
carrying out an analytical study of each competition model. Section 4 is devoted to the presentation of the designed algorithms for the detection
and the refinement of points lying on the separatrix manifold. In Section 5 we describe the Partition of Unity method used for
approximating such curves and surfaces. Section 6 shows some numerical results in both the two and three dimensional phase spaces.

\section{The two populations model}

Let us consider the following competition model, with $N$ and $E$ denoting the red native and the grey exotic squirrels, respectively,

\begin{equation} \label{model}
\begin{array}{ll}
\frac{ \displaystyle dN}{ \displaystyle dt}  =p \bigg(1- \frac{ \displaystyle N}{ \displaystyle u} \bigg)N-aE(1-b)N,  & \textrm{} \\
\vspace{.01cm}\\
\frac{  \displaystyle dE}{  \displaystyle dt}=r \bigg(1- \frac{ \displaystyle E}{ \displaystyle z} \bigg)E-cN(1-b)E,  & \textrm{} \\
\end {array} 
\end{equation}
where $p$ and $r$ are the growth rates of $N$ and $E$, respectively, $a$ and $c$ are their competition rates,
$u$ and $z$ are the respective carrying capacities of the two populations and $b$ denotes the fraction of red squirrels which hide
in a niche, unreachable by the exotic population. Model (\ref{model}) describes the interaction of the two different populations of squirrels within the
same environment.

There are four equilibria for the model, which are given by the points
\begin{displaymath}
\left.
\begin{array}{ll}
\vspace{0.1cm}
E_0=(0,0); \hskip 0.2cm E_1=(0,z);\hskip 0.2cm E_2=(u,0); \\
E_3= \bigg (\displaystyle \frac{ ur[p-az(1-b)]}{pr-aucz(1-b)^2}, \frac{zp[r-cu(1-b)]}{pr-aucz(1-b)^2} \bigg).
\end{array}
\right.
\end{displaymath}
Apart from the presence of the niche, the model is a classical one, thus
assuming that all parameters are positive, we barely summarize the  stability and feasibility results for the equilibrium points   in Table
\ref{tabella}.
\begin{table}[!htbp]
\begin{center}
\begin{tabular}{c|c|c}
\hline
\rule[0mm]{0mm}{3ex}
Equilibrium &  Feasibility & Stability  \\ 
\hline
\rule[0mm]{-1mm}{3ex}
\bfseries  $E_0$  & always feasible & unstable\\
\bfseries  $E_1$  & always feasible &  $p<az(1-b)$ \\ 
\bfseries  $E_2$  & always feasible & $r<cu(1-b)$ \\
\bfseries  $E_3$  &$ p>az(1-b),$ $ r> cu(1-b),$ & $r>cu(1-b), $ $p>az(1-b)$\\
\bfseries   &or  &\\
\bfseries    & $p<az(1-b),$ $r< cu(1-b)$  & \\
\hline
\end{tabular} 
\end{center}
\vspace{0.3cm}
\caption{Feasibility and stability conditions for the equilibria of the system \eqref{model}.}
\label{tabella}
\end{table}
\\
As we can deduce from Table \ref{tabella}, with the parameter
values: $r=1$, $p=2$, $b=0.5$, $u=1$, $c=3$,  $a=2$, $z=3$, the origin $E_0$ is an unstable equilibrium,
$E_1$ and $E_2$ are stable equilibria, and $E_{3}$ is an unstable equilibrium, specifically a saddle point.
In this situation the competitive exclusion principle applies, i.e. only one population survives.
This suggests the existence of a separating curve that divides the model domain into two subregions, called
basins of attraction of each respective equilibrium, each containing paths tending to either $E_1$ or $E_2$.

\begin{figure}[ht!]
\begin{center}
  \includegraphics[height=.3\textheight]{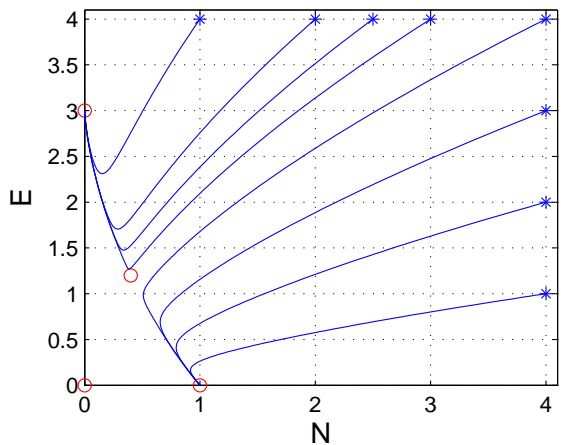}
  \caption{Example of initial conditions and trajectories converging to equilibria for the model (\ref{model}).}
\label{ritratto2d}
\end{center}
\end{figure}

\begin{figure}[ht!]
\begin{center}
  \includegraphics[height=.3\textheight]{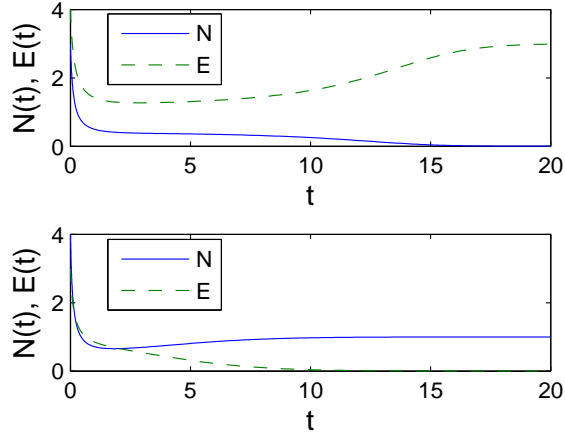}
  \caption{Time evolution of the two populations for the initial conditions $ \boldsymbol{x}_
{3}=(3,4)$ (top) and  $ \boldsymbol{x}_{5}=(4,3)$ (bottom).}
\label{ritratto2d_time}
\end{center}
\end{figure}

In Figure \ref{ritratto2d} we show trajectories starting from the initial conditions $ \boldsymbol{x}_{1}=(1,4)$,
$ \boldsymbol{x}_{2}=(2,4)$, $ \boldsymbol{x}_
{3}=(3,4)$, $ \boldsymbol{x}_{4}=(4,4)$,  $ \boldsymbol{x}_{5}=(4,3)$, $ 
\boldsymbol{x}_{6}=(4,2)$, $ \boldsymbol{x}_{7}=(2.5,4)$ and $ \boldsymbol{x}_{8}=
(4,1)$, and
converging to the point $E_1$ of coordinates $\left(0, 3 \right)$ and to the point $E_2$ of
coordinates $\left( 1,0 \right)$. The time evolution of populations with initial conditions  $ \boldsymbol{x}_{3}=(3,4)$
and  $ \boldsymbol{x}_{5}=(4,3)$ is shown in Figure \ref{ritratto2d_time}, top and bottom respectively,
for the same parameter set.


\section{The three populations model}

Let us now consider the three population competition model, with $N$, $A$ and $E$ denoting the
red native, the red indigenous and the grey exotic squirrels, respectively,

\begin{equation} \label{model3d}
\begin{array}{ll}
\frac{ \displaystyle  dN}{ \displaystyle  dt}=p \bigg (1- \frac{ \displaystyle N}{ \displaystyle u} \bigg )N-aE(1-b)N,  & \textrm{} \\
\vspace{.01cm}\\
\frac{ \displaystyle  dA}{ \displaystyle  dt}=q \bigg (1- \frac{ \displaystyle  A}{ \displaystyle  v} \bigg )A-cE(1-e)A,  & \textrm{} \\
\vspace{.01cm}\\
\frac{ \displaystyle  dE}{ \displaystyle dt}=r \bigg (1- \frac{ \displaystyle  E}{ \displaystyle  z} \bigg)E-fN(1-b)E-gA(1-e)E,  & \textrm{} 
\end{array}
\end{equation}
where $p$, $q$ and $r$ are the growth rates of $N$, $A$ and $E$, respectively, $a$, $c$, $f$ and $g$ are the competition rates,
$u$, $v$ and $z$ are the carrying capacities of the three populations, $b$ and $e$ denote the fraction of the populations
$N$ and $A$, respectively, which hide in a niche. 
Model (\ref{model3d}) describes the interaction of the three different
populations of squirrels, in which there is no competition between the two red squirrel populations because they are
assumed to occupy different habitats \cite{grosso12}.

The analytical study of the model show that the critical points are given by
\begin{displaymath}
\left.
\begin{array}{ll}
\vspace{0.1cm}
E_0 = &(0,0,0); \hskip 0.2cm E_1 = (u,0,0); \hskip 0.2cm E_2 = (0,v,0); \hskip 0.2cm E_3 = (u,v,0); \hskip 0.2cm E_4 = (0,0,z);\\
\vspace{0.1cm}
E_5 = &\bigg( \frac{\displaystyle {ur[az(1-b)-p]}}{\displaystyle{azuf(1-b)^2-pr}},0,\frac{\displaystyle{zp[fu(1-b)-r]}}{\displaystyle{azuf(1-b)^2-pr}}\bigg);\\
\vspace{0.1cm}
E_6 =&\bigg(  0,\frac{\displaystyle{vr[cz(1-e)-q]}}{\displaystyle{czvg(1-e)^2-qr}},\frac{\displaystyle{zq[vg(1-e)-r]}}{\displaystyle{czvg(1-e)^2-qr}}\bigg);\\
\vspace{0.1cm}
E_7 = & \bigg( \frac{\displaystyle{u[ \alpha-prq-zaq(1-b)(vg(1-e)-r)]}}{\displaystyle{\alpha+\beta-prq}},\\
\vspace{0.1cm}
& \frac{\displaystyle{v[ \beta-prq-pcz(1-e)(fu(1-b)-r)]}}{\displaystyle{\alpha+\beta-prq}},\\
& \frac{\displaystyle{zpq[fu(1-b)+vg(1-e)]}}{\displaystyle{\alpha+\beta-prq}}\bigg).
\end{array}
\right.
\end{displaymath}
where, for brevity, we indicated for $E_7$, 
\begin{displaymath}
\alpha = pczvg(1-e)^2, \hskip 0.2cm \beta = azufq(1-b)^2.
\end{displaymath}
Like in the 2D case, we assume that all parameters are positive and we summarize the  study of   stability and feasibility of the equilibrium points in Table \ref{tabella1}.
\begin{table}[!htbp]
\begin{center}
\begin{tabular}{c|c|c}
\hline
\rule[0mm]{0mm}{3ex}
Equilibrium &  Feasibility & Stability  \\ 
\hline
\rule[0mm]{-1mm}{3ex}
\bfseries $E_0$  & always feasible &   unstable\\
\bfseries  $E_1$  & always feasible &   unstable\\ 
\bfseries  $E_2$  & always feasible &   unstable\\
\bfseries  $E_3$  & always feasible &  $r<fu(1-b) + vg(1-e) $ \\
\bfseries  $E_4$  & always feasible &   $q<cz(1-e)$, $  p<az(1-b)$\\
\bfseries  $E_5$  & $r>fu(1-b)$, &  $r>fu(1-b)$, $ p>az(1-b),$\\
\bfseries &  $p>az(1-b),$ &  $c(1-e)zp[fu(1-b)-r] < q[azuf(1-b)^2-pr]$\\
\bfseries &or & \\
\bfseries & $r<fu(1-b),$ &  \\
\bfseries & $p<az(1-b)$ & \\
\bfseries  $E_6$  &$q>cz(1-e)$,  &  $q>cz(1-e)$, $ r>vg(1-e),$\\
\bfseries & $ r>vg(1-e),$ &   $ a(1-b)zq[vg(1-e)-r] < p[czvg(1-e)^2-qr]$ \\
\bfseries & or &\\
\bfseries &$q<cz(1-e)$, & \\
\bfseries & $ r<vg(1-e)$ & \\
\hline
\end{tabular}
\end{center}
\vspace{0.3cm}
\caption{Feasibility and stability conditions for the equilibria of the system \eqref{model3d}.}
\label{tabella1}
\end{table}
As we can deduce from this table with the choice of the parameters $r = 9$, $q = 0.6$, $p = 0.6$, $b = 0.5$,
$u = 1.5$, $c = 8$, $a = 8$, $z = 3$, $v = 2$, $e = 0.5$, $f = 6$, $g = 5$, the points $E_3$ and $E_4$ are stable
equilibria, $E_5$ and $E_6$ are not feasible. We verify numerically that with this choice $E_7$ is a saddle point.
This suggests the existence of a separating surface that divides the model domain into two basins of
attraction, each of them containing one path tending to $E_3$ or $E_4$.

\begin{figure}[ht!]
\begin{center}
  \includegraphics[height=.3\textheight]{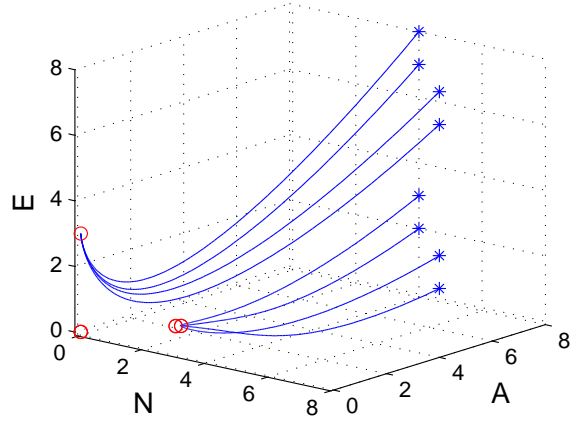}
  \caption{Example of initial conditions and trajectories converging to equilibria for the model
 problem (\ref{model3d}).}
\label{ritratto_fase_3D}
\end{center}
\end{figure}

In Figure \ref{ritratto_fase_3D} we show trajectories starting from the initial conditions $\boldsymbol{x}_{1}=(4,8,3)$, $\boldsymbol{x}_{2}=(4,8,2)$, $\boldsymbol{x}_
{3}=(4,8,7)$, $\boldsymbol{x}_{4}=(4,8,8)$, $\boldsymbol{x}_{5}=(8,4,3)$, $\boldsymbol{x}_{6}=(8,4,2)$, $\boldsymbol{x}_{7}=(8,4,7)$ and $\boldsymbol{x}_{8}
=(8,4,8)$, and converging to the point $E_3$ of coordinates $\left(1.5, 2, 0 \right)$ and to the point $E_4$ of coordinates $\left( 0,0, 3\right)$. The time evolution of populations with initial conditions  $ \boldsymbol{x}_
{2}=(4,8,2)$ and  $ \boldsymbol{x}_{5}=(8,4,8)$ is shown in Figure \ref{ritratto_fase_3D}, top and bottom respectively.

\begin{figure}[ht!]
\begin{center}
\includegraphics[height=.3\textheight]{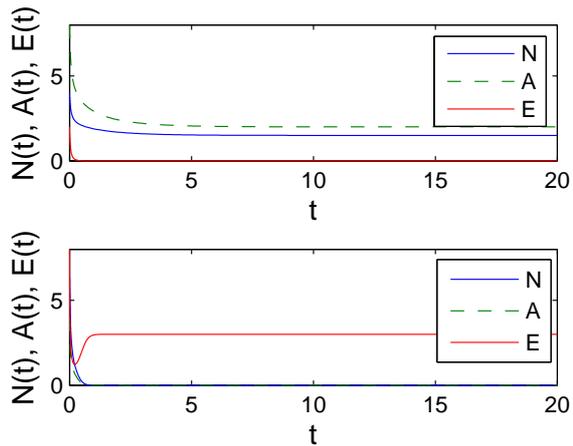}
  \caption{Time evolution of the three populations for the initial conditions $ \boldsymbol{x}_{2}=(4,8,2)$ (top) and  $ \boldsymbol{x}_{5}=(8,4,8)$ (bottom).}
\label{ritratto_fase_3D_time}
\end{center}
\end{figure}


\section{Detection and refinement of separatrix points}

At first, to determine the separatrix curve and surface for (\ref{model}) and (\ref{model3d}), respectively, we need to consider a set of points as initial
conditions in a square domain $[0,\gamma]^2$, where $\gamma \in \RR^+$, and in a cube domain $[0,\gamma]^3$.
Then we take points in pairs and we check if trajectories of the two points converge to different
equilibria. If this is the case, then we proceed with a bisection algorithm to determine a separatrix
point.
Then we perform a refinement of the set of separatrix points. In fact, in general we find a large number of separatrix points. In what follows
we propose a refinement
process which computes a smaller set of points. The set of the refined points is then interpolated using a suitable method
described in Section 5.

Let now qD denote the dimension of the phase space of the dynamical system.
In the 2D case we start considering $n$ equispaced initial conditions on each edge of the
square $[0,\gamma]^2$  and the bisection algorithm is applied with the following initial conditions 
$$
(x_i,0) \quad \textrm{and} \quad (x_i,\gamma), \quad i=1, \ldots, n, 
$$
$$
(0,y_i) \quad \textrm{and} \quad (\gamma,y_i), \quad i=1, \ldots, n. 
$$
Performing the bisection algorithm, a certain number of points, in general large, is found on the separatrix curve. 
The $N$ points found by the bisection algorithm are collected in a matrix 
$A=(a_{j,k})$, $j=1, \ldots, N$, $k=1,2$, and then refined in order to obtain a smaller set of well distributed nodes on the separatrix curve. 
So we divide the interval $[0,M]$, where $M=  \max_{ \displaystyle j} (a_{j,1}), $ $ j=1, \ldots, N,$ in $L$ subintervals and we compute an average of the found points in each subinterval. Given a vector of equispaced points $ x_l$, $l=1, \ldots, L+1$, in the interval $[0,M]$, we define
$$
I_{l}= \{j : a_{j,1} \in [x_l,x_{l+1}) \},
$$
with $l=1, \ldots, L$.
Starting from the matrix $A=(a_{j,k})$, $ j=1, \ldots, N $, $k=1,2$, we define the matrix of the refined points $A^{'}=(a^{'}_{j,k})$, whose entries are given by 
\begin{displaymath}
a^{'}_{1,k}=a_{1,k}, \quad k=1,2,
\end{displaymath}
\begin{displaymath}
a^{'}_{j,1}= \frac{\sum_{j \in I_{l}}  a_{j,1} }{ Card(I_{l})}, \quad l=1, \ldots ,L,
\end{displaymath}

\begin{displaymath}
a^{'}_{j,2}= \frac{ \sum_{j \in I_{l}}  a_{j,2} }{ Card(I_{l})}, \quad l=1, \ldots ,L,
\end{displaymath}

\begin{displaymath}
a^{'}_{K+2,k}=a_{N,k}, \quad k=1,2,
\end{displaymath}
and $j=2, \ldots,K+1$, where $K$ is the number of subintervals containing at least a point.
We summarize the steps in Algorithm 1.

\begin{table}[!ht]
\begin{center}
\begin{tabular}{p{13cm}*{1}{c}}
\hline
\rule[0mm]{0mm}{3ex}
\textsl{Algorithm 1.}  \\ 
\hline
\\[\smallskipamount]
{\fontfamily{pcr} \selectfont Step 1:}
 Definition of initial conditions.
\vskip 0.06 cm
\hskip 2.2 cm (Equispaced vectors on edges of the square). 
\vskip 0.06 cm
\hskip 2.2 cm $P_{i}^{1}=(x_i,0)$, $P_{i}^{2}=(x_i,\gamma)$, $P_{i}^{3}=(0,y_i)$, $P_{i}^{4}=(\gamma,y_i)$,
\vskip 0.06 cm
\hskip 2.2 cm $i=1, \ldots, n$.
\vskip 0.06 cm
\hskip 2.2 cm Set $j=0$,
\vskip 0.06 cm
\hskip 2.9 cm $s=1$.
\vskip 0.2 cm
 {\fontfamily{pcr} \selectfont Step 2:} While $s<=3$
\vskip 0.2 cm
\hskip 0.55 cm
{\fontfamily{pcr} \selectfont Step 3:} For $i=1, \ldots, n$
\vskip 0.2 cm
\hskip 1.2 cm
\noindent {\fontfamily{pcr} \selectfont Step 4:}  
If $P_{i}^{s} \rightarrow E_1$ and $P_{i}^{s+1} \rightarrow E_2$ or vice versa 
\vskip 0.02 cm
\hskip 3.4 cm then $j=j+1$,
\vskip 0.2 cm
\hskip 4.1 cm $Q_{j,k}=BISECTION(P_{i}^s,P_{i}^{s+1})$, $k=1,2$.
\vskip 0.2 cm
\hskip 2.2 cm
 $s=s+2$.
\vskip 0.2 cm
{\fontfamily{pcr} \selectfont Step 5:} 
Define $N = j$
\vskip 0.06 cm
\hskip 2.9 cm  (number of points found by the bisection algorithm),
\vskip 0.06 cm 
\hskip 2.9 cm $M= \max_{j} Q_{j,1}$, $j=1, \ldots, N$.
\vskip 0.06 cm
\hskip 2.2 cm Define an equispaced vector $x_l$, $l=1, \ldots , L+1$.
\vskip 0.06 cm
\hskip 2.2 cm Set $I_{l}= \{j : Q_{j,1} \in [x_l,x_{l+1}) \}$, $l=1, \ldots , L$.
\vskip 0.06 cm
\hskip 2.2 cm $QQ_{1,k}=Q_{1,k}$, \quad  $k=1,2$,
\vskip 0.06 cm
\hskip 2.2 cm $QQ_{j,k}= \displaystyle \frac{  \sum_{j \in I_{l}}Q_{j,k} }{  Card(I_{l})},  \quad l=1, \ldots , L,   \quad k=1,2$,
\vskip 0.06 cm
\hskip 2.2 cm $QQ_{K+2,k}=Q_{N,k}$, \quad $k=1,2$,
\vskip 0.06 cm
\hskip 2.2 cm and $j=2, \ldots,K+1$, where $K$ is the number of subintervals
\vskip 0.06 cm
\hskip 2.2 cm   containing at least a point.
\vskip 0.2 cm
{\fontfamily{pcr} \selectfont Step 6:} 
INTERPOLATION($QQ_{j,k}^1$), where $QQ_{j,k}^1$ is the set composed by
\vskip 0.06 cm
\hskip 2.2 cm  the points  found by the refinement algorithm, saddle and origin.
\vskip 0.06 cm
\hskip 2.2 cm (See Section 5).\\
\\[\smallskipamount]
\hline
\end{tabular}
\label{algo1}
\end{center}
\vspace{0.3cm}
\end{table}

For the 3D case, we use a similar technique. At first we construct a grid on the faces of the cube and the bisection algorithm is applied with the following initial conditions 
$$
(x_{i_1},y_{i_2},0) \quad \textrm{and} \quad (x_{i_1},y_{i_2},\gamma), \quad i_1=1, \ldots, n, \quad i_2=1, \ldots, n,
$$
$$
(x_{i_1},0,z_{i_2}) \quad \textrm{and} \quad (x_{i_1},\gamma,z_{i_2}), \quad i_1=1, \ldots, n, \quad i_2=1, \ldots, n,
$$
$$
(0,y_{i_1},z_{i_2}) \quad \textrm{and} \quad (\gamma,y_{i_1},z_{i_2}), \quad i_1=1, \ldots, n, \quad i_2=1, \ldots, n.
$$
The $N$ points found by the bisection algorithm are organized in a matrix 
$A=(a_{j,k})$, $j=1, \ldots, N$, $k=1,2,3$, and then refined. 
We define
$$
M_x= \max_{ \displaystyle j} (a_{j,1}), \quad j=1, \ldots, N,
$$
$$
M_y= \max_{ \displaystyle j} (a_{j,2}), \quad j=1, \ldots, N,
$$\\
and we divide the interval $[0,M_x]$ in $L$ subintervals  and $[0,M_y]$ in $H$ and we make an average of the points in each subinterval. 
Given a vector of equispaced points $x_l$, $l=1, \ldots , L+1$, in $[0,M_x]$, and a vector $y_h$, $h=1, \ldots, H+1$, in $[0,M_y]$, let us define 
$$
I_{lh}= \{j : a_{j,1} \in [x_l,x_{l+1}] \quad \textrm{and} \quad a_{j,2} \in  [y_{h},y_{h+1}] \},
$$\\
with $l=1, \ldots , L,$ $ h=1, \ldots, H.$ Starting from the matrix $A=(a_{j,k})$ we find the matrix of the refined points $A^{'}=(a^{'}_{j,k})$, whose entries are given by 
$$
a^{'}_{j,1}= \frac{ \sum_{j \in I_{lh}} a_{j,1} }{ Card(I_{lh})},  \quad l=1, \ldots , L, \quad h=1, \ldots ,  H,
$$
$$
a^{'}_{j,2}= \frac{ \sum_{j \in I_{lh}} a_{j,2} }{ Card(I_{lh})}, \quad l=1, \ldots , L, \quad h=1, \ldots ,  H,
$$
$$
a^{'}_{j,3}= \frac{ \sum_{j \in I_{lh}} a_{j,3} }{ Card(I_{lh})}, \quad l=1, \ldots , L, \quad h=1, \ldots ,  H,
$$
and $j=1, \ldots,K$, where $K$ is the number of subintervals containing at least a point.
We summarize the steps in Algorithm 2.

\vspace{0.1cm}
\begin{table}[!ht]
\begin{center}
\begin{tabular}{p{13.1cm}*{1}{c}}
\hline
\rule[0mm]{0mm}{3ex}
\textsl{Algorithm 2.}  \\ 
\hline
\\[\smallskipamount]
{\fontfamily{pcr} \selectfont Step 1:}
 Definition of initial conditions.
\vskip 0.06 cm
\hskip 2.2 cm  (Grid on the faces of the cube). 
\vskip 0.06 cm
\hskip 2.2 cm $P_{i_1,i_2}^{1}=(x_{i_1},y_{i_2},0)$, $P_{{i_1},{i_2}}^{2}=(x_{i_1},y_{i_2},\gamma)$, $P_{{i_1},{i_2}}^{3}=(x_{i_1},0,z_{i_2})$,
\vskip 0.06 cm
\hskip 2.2 cm $P_{{i_1},{i_2}}^{4}=(x_{i_1},\gamma,z_{i_2})$, $P_{{i_1},{i_2}}^{5}=(0,y_{i_1},z_{i_2})$, $P_{{i_1},{i_2}}^{6}=(\gamma,y_{i_1},z_{i_2})$,
\vskip 0.06 cm
\hskip 2.2 cm ${i_1}=1, \ldots, n$, \hskip 0.2 cm ${i_2}=1, \ldots, n.$
\vskip 0.06 cm
\hskip 2.2 cm Set $j=0$,
\vskip 0.06 cm
\hskip 2.9 cm $s=1$.
\vskip 0.2 cm
{\fontfamily{pcr} \selectfont Step 2:} While $s<=5$
\vskip 0.2 cm
\hskip 0.55 cm {\fontfamily{pcr} \selectfont Step 3:} For $i_1=1, \ldots, n$
\vskip 0.2 cm
\hskip 1.2 cm {\fontfamily{pcr} \selectfont Step 4:} For $i_2=1, \ldots, n$
\vskip 0.2 cm
\hskip 1.9 cm {\fontfamily{pcr} \selectfont Step 5:}  
 If $P_{{i_1},{i_2}}^{s} \rightarrow E_3$ and $P_{{i_1},{i_2}}^{s+1} \rightarrow E_4$ or vice versa
\vskip 0.02 cm
\hskip 4.1 cm then $j=j+1$,
\vskip 0.2 cm
\hskip 4.3 cm $Q_{j,k}=BISECTION(P_{{i_1},{i_2}}^s,P_{{i_1},{i_2}}^{s+1})$, $k=1,2,3$.
\vskip 0.2 cm
\hskip 3.2 cm 
$s=s+2$.
\vskip 0.2 cm
{\fontfamily{pcr} \selectfont Step 6:} 
Define $N = j$
\vskip 0.06 cm 
\hskip 2.95 cm 
(number of points found by the bisection algorithm),
\vskip 0.06 cm 
\hskip 2.95 cm $M_x= \max_{j} Q_{j,1}$, $j=1, \ldots, N$, 
\vskip 0.06 cm 
\hskip 2.95 cm $M_y= \max_{j} Q_{j,2}$,  $j=1, \ldots, N$.
\vskip 0.06 cm 
\hskip 2.2 cm Define two equispaced vectors $x_l$, $l=1, \ldots , L+1$ and $y_h$, 
\vskip 0.06 cm
\hskip 2.2 cm $h=1, \ldots , H+1$.
\vskip 0.06 cm
\hskip 2.2 cm Define $I_{lh}= \{j : Q_{j,1} \in [x_l,x_{l+1}] \quad \textrm{and} \quad Q_{j,2} \in  [y_h,y_{h+1}] \}$, 
\vskip 0.06 cm
\hskip 2.2 cm $l=1, \ldots , L$,  $h=1, \ldots, H$.
\vskip 0.06 cm
\hskip 2.2 cm $QQ_{j,k}= \displaystyle  \frac{ \sum_{j \in I_{lh}}Q_{j,k} }{ Card(I_{lh})},  \quad l=1,\ldots,L, \hskip 3.mm h=1,\ldots,H, \hskip 3.mm k=1,2,3,$
\vskip 0.06 cm
\hskip 2.2 cm and $j=1, \ldots,K$, where $K$ is the number of subintervals 
\vskip 0.06 cm
\hskip 2.2 cm containing  at least a point.
\vskip 0.06 cm
\vskip 0.2 cm
{\fontfamily{pcr} \selectfont Step 7:}
INTERPOLATION($QQ_{j,k}^1$), where $QQ_{j,k}^1$ is the set composed by the 
\vskip 0.06 cm
\hskip 2.2 cm  points found by the refinement algorithm, saddle and origin. (See 
\vskip 0.06 cm
\hskip 2.2 cm Section 5).\\
\\[\smallskipamount]
\hline
\end{tabular}
\label{algo2}
\end{center}
\vspace{0.3cm}
\end{table}


\section{Reconstruction of separatrix curves and surfaces} 

In this section we present the interpolation method, namely the so-called Partition of Unity method, we use to connect the found points applying the refinement algorithm.

Let us consider a set ${\cal X}=\{ \textbf{x}_i, i=1,\ldots,n\}$ of distinct data points
arbitrarily distributed on $\Omega \subseteq \RR^m$, and an associated
set ${\cal F}=\{ f_i, i=1,\ldots,n\}$ of data values.

The basic idea of the Partition of Unity method is to start with a
partition of the open and bounded domain $\Omega \subseteq \RR^m$ into $d$
 cells (subdomains) $\Omega_j$ such that $\Omega \subseteq
\bigcup_{j=1}^{d} \Omega_j$ with some mild overlap among the cells. At
first, we choose a partition of unity, i.e. a family of compactly
supported, non-negative, continuous functions $W_j$ with $\text{supp}(W_j)
\subseteq \Omega_j$ such that 
\begin{equation}
\sum_{j=1}^{d} W_j(\textbf{x}) = 1, \hspace{1cm} \textbf{x} \in \Omega. \nonumber
\end{equation}
Then, for each cell $\Omega_j$ we consider a local approximant $R_j$
and form the global approximant given by
\begin{eqnarray}
\label{pui}
	{\cal I}(\textbf{x})= \sum_{j=1}^{d} R_j(\textbf{x}) W_j(\textbf{x}), \hspace{1cm} \textbf{x} \in \Omega.
\end{eqnarray}
Note that if the local approximants satisfy the interpolation conditions
at data point $\textbf{x}_i$, i.e.  $R_j(\textbf{x}_i)=f_i$, then the global approximant
also interpolates at this node, i.e. ${\cal I}(\textbf{x}_i)=f_i$, for
$i=1,\ldots,n$ (see \cite{Fasshauer07,Wendland05} for further details).

As a local approximant we can take a radial basis function interpolant $R_j:\Omega
\rightarrow \RR$, which has the form
\begin{eqnarray}
\label{intfun}
R_j({\textbf{x}})=\sum_{j=1}^n \alpha_j \phi (||\textbf{x}-\textbf{x}_j||_2). \hskip0.5cm \textbf{x}\in \Omega,
\end{eqnarray}
Here $\phi:[0,\infty) \rightarrow \RR$ is called a \textsl{radial basis function}, $||\cdot||_2$ is
the Euclidean norm, and $\{\alpha_j\}$ are the coefficients to be determined by solving the linear
system generated by radial basis functions. Moreover, $R_j$ satisfies the interpolation conditions
$R_j({\textbf{x}}_i)=f_i$, $i=1,\ldots,n$ (see \cite{Buhmann03,Fasshauer07}).

Usually, it can be highly advantageous to work with locally supported functions
since they lead to sparse linear systems. In \cite{Wendland05}) Wendland found
a class of radial basis functions which are smooth, locally supported, and
strictly positive definite on $\RR$. They consist of a
product of a truncated power function and a low degree polynomial. For
example, here we take the Wendland $C2$ function
\begin{equation}
\phi(r)  =  \displaystyle{\left(1-\beta r\right)_+^4\left(4\beta r+1\right)},
\nonumber
\end{equation}
where $r=||\textbf{x}-\textbf{x}_j||_2$, $\beta \in \RR^+$ is the shape parameter, and $(\cdot)_+$ denotes the truncated
power function. This means that the function $\phi(r)$ is nonnegative; in fact, $(1-\beta r)_+$ is defined as $(1-\beta r)$
for $r\in \left[0,1/\beta\right]$, and 0 for $r > 1/\beta$.

In order to be able to formulate error bounds we need some technical conditions. Then, we require the partition of unity functions $W_j$ to be \textsl{k-stable}, i.e. each $W_j \in C^k(\RR^m)$ and for every multi-index $\mu \in \NN_0^m$ with $|\mu| \leq k$ there exists a constant $C_{\mu} > 0$ such that
\begin{equation}
	\left\|D^{\mu}W_j\right\|_{L_{\infty}(\Omega_j)}\leq C_{\mu}/\delta_j^{|\mu|}. \nonumber
\end{equation}
where $\delta_j$ = diam($\Omega_j$).

Now, after defining the space $C_{\nu}^k(\RR^m)$ of all functions $f \in C^k$ whose derivatives of order $|\mu|=k$ satisfy $D^{\mu}f(\textbf{x})= O(||\textbf{x}||_2^{\nu})$ for $||\textbf{x}||_2 \rightarrow 0$, we have the following approximation theorem \cite{Wendland02}.
\begin{theorem}
	Let $\Omega \subseteq  \RR^m$ be open and bounded and suppose that ${\cal X} = \{\textbf{x}_i, i=1,$ $\ldots,n \}\subseteq \Omega$. Let $\phi \in C_{\nu}^k(\RR^m)$ be a strictly positive definite function. Let $\{\Omega_j\}_{j=1}^{d}$ be a regular covering for $(\Omega, {\cal X})$ and let $\{W_j\}_{j=1}^{d}$ be $k$-stable for $\{\Omega_j\}_{j=1}^{d}$. Then the error beetween $f \in {\cal N}_{\phi}(\Omega)$, where ${\cal N}_{\phi}$ is the native space of $\phi$, and its partition of unity interpolant (\ref{pui}) can be bounded by
\begin{equation}
	|D^{\mu}f(\textbf{\textbf{x}}) - D^{\mu}{\cal I}(\textbf{\textbf{x}})| \leq C h_{{\cal X}, \Omega}^{(k+\nu)/2 - |\mu|} |f|_{{\cal N}_{\phi}(\Omega)}, \hspace{1.cm} \forall \textbf{x} \in \Omega, \quad |\mu| \leq k/2, \nonumber
\end{equation}
	$h_{{\cal X}, \Omega}$ being the so-called \textsl{fill distance}, whose definition is given by $h_{{\cal X}, \Omega} = \sup_{\textbf{x} \in \Omega}$ $\min_{\textbf{x}_j\in {\cal X}} ||\textbf{x}-\textbf{x}_j||_2$. 
\end{theorem}

Some additional assumptions on regularity of $\Omega_j$ are required:
\begin{itemize}
	\item for each $\textbf{x} \in \Omega$ the number of subdomains $\Omega_j$ with $\textbf{x} \in \Omega_j$ is bounded by a global constant $K$;
	\item each subdomain $\Omega_j$ satisfies an interior cone condition (see \cite{Wendland05});
	\item the local fill distances $h_{{\cal X}_j, \Omega_j}$ are uniformly bounded by the global fill distance $h_{{\cal X}, \Omega}$, where ${\cal X}_j={\cal X} \cap \Omega_j$.
\end{itemize}

Note that the Partition of Unity method preserves the local approximation order for the global fit.
Hence, we can efficiently compute large radial basis function
interpolants by solving small radial basis functions interpolation problems (in parallel as well)
and then combine them together with the global partition of unity
$\{W_j\}_{j=1}^{d}$. This approach enables us to decompose a large problem into many small problems,
and at the same time ensures that the accuracy obtained for the local fits is carried over to the
global fit. In particular, the Partition of Unity method can be thought as a Shepard's method with higher-order data, since
local approximations $R_j$ instead of data values $f_j$ are used. Moreover, the use of Wendland's
functions guarantees a good compromise between accuracy and stability.


\section{Numerical experiments} 

\begin{figure}[ht!]
\begin{center}
  \includegraphics[height=.22\textheight]{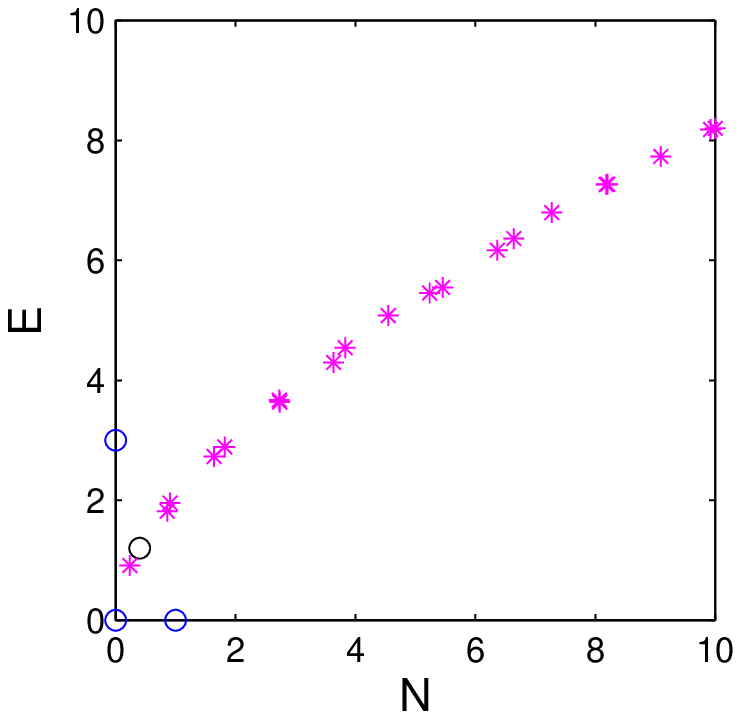} 
  \includegraphics[height=.22\textheight]{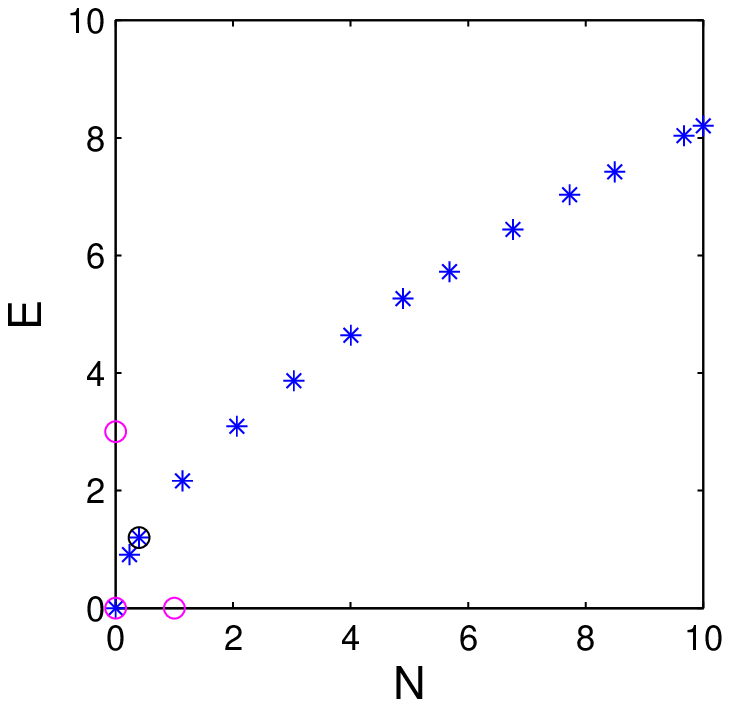} 
    \caption{Set of points detected by the bisection algorithm (left) and set of points found by the refinement algorithm (right) in the 2D case. The circles represent equilibria and saddle points.}
\label{figura1}
\end{center}
\end{figure}

\begin{figure}[ht!]
\begin{center}
  \includegraphics[height=.3\textheight]{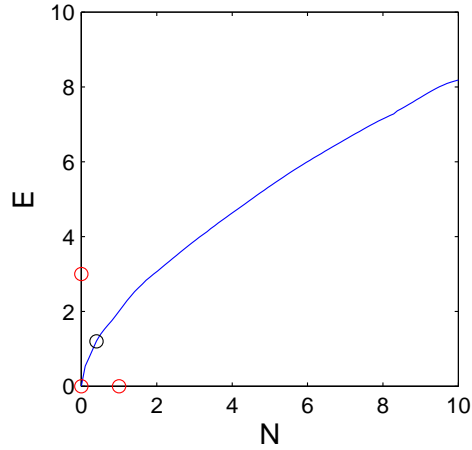} 
  \caption{Reconstruction of the separatrix curve using $\beta=0.025$.}
\label{figura2}
\end{center}
\end{figure}

In this section we summarize the extensive experiments to test our detection and approximation techniques.  In the following we fix $\gamma=10$.  We first
refer to the dynamical system (\ref{model}), taking $r=1$, $p=2$, $b=0.5$, $u=1$, $c=3$,  $a=2$, and $z=3$.  
For example, Figure \ref{figura1} (left) shows the points found using $n=12$. Dividing the interval $[0,M]=[0,10]$ in $L=10$ subintervals and considering the $N=20$ points picked up on the separatrix curve, the refinement process provides us the $K+2=L+2=12$ points.
To this set we add the origin and the saddle point, as shown in Figure \ref{figura1} (right).
A crucial task for the accuracy of the Partition of Unity method concerns
the choice of the shape parameter $\beta$ of Wendland's function. In fact, it can significantly affect the approximation result and, therefore, the quality of the separatrix curve. From our study we found that good shape parameter values are given for $0.001 \le \beta \le 0.04 $. In Figure \ref{figura2}  we show the curve obtained by approximating the refined data set when we consider the value $\beta=0.025$ as shape parameter for the Wendland $C2$ function and a number $d=3$ of partitions of $\Omega$.\\

\begin{figure}[ht!]
\begin{center}
  \includegraphics[height=.22\textheight]{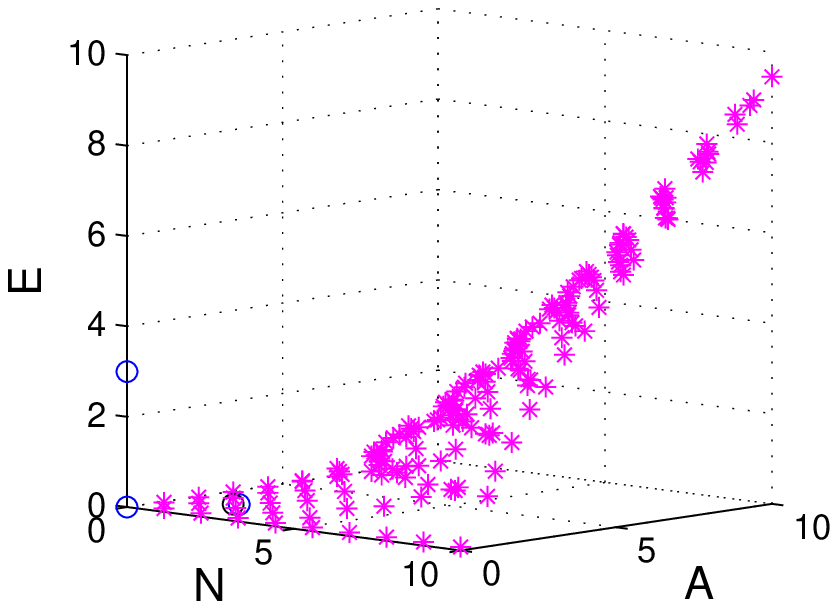}
  \includegraphics[height=.22\textheight]{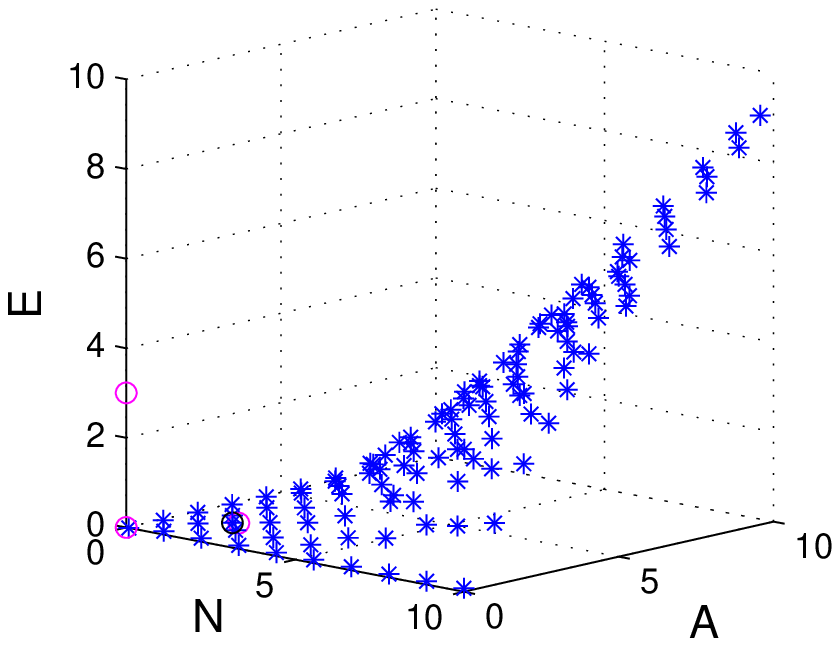} 
    \caption{Set of points detected by the bisection algorithm (left) and set of points found by the refinement algorithm (right) in the 3D case. The circles represent equilibria and saddle points.}
\label{figura3}
\end{center}
\end{figure}

\begin{figure}[ht!]
\begin{center}
  \includegraphics[height=.3\textheight]{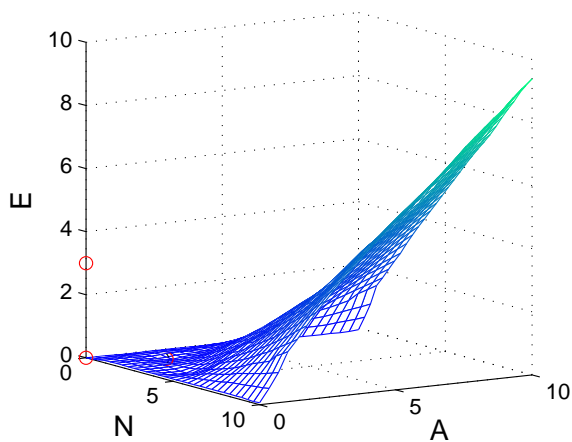} 
  \caption{Reconstruction of the separatrix surface using $\beta=0.005$.}
\label{figura4}
\end{center}
\end{figure}

Then, we consider the dynamical system (\ref{model3d}), taking  $r = 9$, $q = 0.6$, $p = 0.6$, $b = 0.5$, $u = 1.5$, $c = 8$, $a = 8$, $z = 3$, $v = 2$, $e = 0.5$, $f = 6$, $g = 5$ as values of biological parameters.
As an example, in Figure  \ref{figura3} (left) we show the points found
choosing $n=10$. The $N=182$ points have been refined taking $H=13$ and
$L=13$. In this way we obtain $K=117$ points. To this set we add the origin
and the saddle point, as shown in Figure \ref{figura3} (right). From our study
we found that good shape parameter values are given for $0.001 \le \beta \le 0.03$.
In Figure \ref{figura4}  we show the surface obtained by approximating the refined
data set when we consider the value $\beta=0.005$ as shape parameter for the Wendland $C2$ function and a number $d=4$ of partitions of $\Omega$.\\

\end{document}